\documentclass[12pt]{amsart}

\usepackage[mathscr]{eucal}
\usepackage{amsmath,amssymb,amscd,amsthm}
\usepackage{color}
\usepackage{amsmath,amsthm,amssymb,graphicx}
\usepackage{epsfig}
\usepackage{graphicx}
\usepackage{tikz}
\graphicspath{ {images/} }
\newtheorem{theorem}{Theorem}

\theoremstyle{definition}

\theoremstyle{plane}

\def \R{\mathbb R}
\def \beq{ \begin{equation} }
\def \eeq{\end{equation}}

\usepackage[top=4cm, bottom=4cm, left=3.5cm, right=3.5cm]{geometry}

\title{Relative equilibria in curved restricted 4-body problems}

\begin{document}
\maketitle

\markboth{Florin Diacu and Sawsan Alhowaity}{Relative equilibria in curved restricted 4-body problems}

\author{\begin{center}
Florin Diacu$^{1,2,3}$ and Sawsan Alhowaity$^{3,4}$\\
\bigskip
$^1$Yangtze Center of Mathematics, Sichuan University, Chengdu, China\\
$^2$Yale-NUS College, National University of Singapore, Singapore\\
$^3$Department of Mathematics and Statistics, 
University of Victoria, Canada\\
$^4$Department of Mathematics, University of Shaqra, Saudi Arabia\\
\medskip
florin.diacu@yale-nus.edu.sg,\ \ salhowaity@su.edu.sa
\end{center}

\begin{abstract}
We consider the curved 4-body problems on spheres and hyperbolic spheres. After obtaining a criterion for the existence of quadrilateral configurations on the equator of the sphere, we study two restricted 4-body problems, one in which two masses are negligible, and another in which only one mass is negligible.
In the former we prove the evidence square-like relative equilibria, whereas in the latter we discuss the existence of kite-shaped relative equilibria. 
\end{abstract}

\section{Introduction}

The classical $N$-body problem has a long history. Isaac Newton first proposed it in 1687 in his first edition of {\it Principia} in the context of the Moon's motion. He assumed that universal gravitation acts between celestial bodies (reduced to point masses) in direct proportion with the product of the masses and in inverse proportion with the square of the distance. The study of the $N$-body problem was further advanced by the Bernoulis, Lagrange, Laplace, Euler, Cauchy, Jacobi, Dirichlet, Poincar\'e and many others.

The idea of extending the gravitational force between point masses to spaces of constant curvature occurred soon after the discovery of hyperbolic geometry. In the 1830s, independently of each other, Bolyai and Lobachevsky realized that there must be an intimate connection between the laws of physics and the geometry of the universe, \cite{Bolyai}, \cite{Lobachevsky}, \cite{Kragh}. A few years earlier, Gauss had interpreted Newton's gravitational  law as stating that the attracting force between bodies is inversely proportional with the area of the sphere of radius equal to the distance between the point masses (i.e.\ proportional to $1/r^2$, where $r$ is the distance). Using this idea, Bolyai and Lobachevsky suggested that, should space be hyperbolic, the attracting force between bodies must be inversely proportional to the hyperbolic area of the corresponding hyperbolic sphere (i.e.\ proportional to $1/\sinh^2(|\kappa|^{1/2}r)$, where $r$ is the distance and $\kappa<0$ the curvature of the hyperbolic space).
This is equivalent to saying that, in hyperbolic space, the potential that describes the gravitational force is proportional to $\coth(|\kappa|^{1/2}r)$.

The above analytic expression of the potential was first introduced by Schering, \cite{Schering1}, \cite{Schering2}, 
and then extended to elliptic space by Killing, [22--24]. 
But with no physical ways of checking the validity of this generalization of the gravitational force, it was unclear whether the cotangent potential had any physical meaning, the more so since Lipschitz had proposed a different extension of the law, which turned out to be short lived, \cite{Lipschitz}. The breakthrough came at the dawn of the 20th century when Liebmann made two important discoveries, \cite{Liebmann1}, \cite{Liebmann2}. He showed that two basic properties of the Newtonian potential are also satisfied by the
cotangent potential: (1) in the Kepler problem, which studies the motion of one body around a fixed centre, the potential is a harmonic function (i.e.\ a solution of the Laplace equation in the Euclidean case, but of the Laplace-Beltrami equation in the non-flat case); (2) in both the flat and the non-flat case, all bounded orbits of the Kepler problem are closed, a property discovered by Bertrand for the Newtonian law, \cite{Bertrand}. These similarities
between the flat and the curved problem convinced the scientific community that the cotangent potential was the natural way to express gravity in spaces of constant curvature. 

The curved $N$-body problem became somewhat neglected after 
the birth of general relativity, but was revived after the discretization
of Einstein's equation showed that an $N$-body problem in spaces
of variable curvature is too complicated to be treated with analytical tools. In the 1990s, the Russian school of celestial mechanics considered both the curved Kepler and the curved 2-body problem,
\cite{Kozlov}, \cite{Shchepetilov}. After understanding that, unlike in the Euclidean case, these problems are not equivalent, the latter failing to be integrable, \cite{Shchepetilov}, the 2-body case was intensively studied by several researchers of this school. More recently, the work of Diacu, Santoprete, and P\'erez-Chavela considered the curved $N$-body problem for $N>2$ in a new framework, leading to many interesting results, [3--20], \cite{Perez}. Other researchers developed these ideas further, \cite{Garcia} \cite{Martinez1}, \cite{Martinez2}, [37--40], and the problem is growing in popularity. 

In this short note we prove three results. The first is a criterion for the existence of quadrilateral relative equilibria on the equator of
the sphere. The second shows that if two masses are negligible and
the other two are equal, then square-like relative equilibria exists on spheres, but---surprisingly---not on hyperbolic spheres. The element of surprise arises from the fact that, in the general problem, square-like equilibria exist both on the hyperbolic sphere and on the sphere (except for the case when they are on the equator), \cite{Diacu03}. Finally we prove that if only one mass is negligible and the other three are equal, some kite-shaped relative equilibria exist on spheres, but not on hyperbolic spheres.

\section{Equations of motion}

We consider the motion of four bodies on 2-dimensional 
surfaces of constant curvature $\kappa$, namely spheres
$\mathbb S_\kappa^2$ for $\kappa>0$, the Euclidean plane 
$\R^2$ for $\kappa=0$, and hyperbolic spheres $\mathbb H_\kappa^2$ for $\kappa<0$. We will arrange these 
objects in $\R^3$ such that they all have a common point
at which lie all the north poles of the spheres and the
vertices of the hyperbolic spheres, to all of which the plane
$\R^2$ is tangent. If we fix the origin of a coordinate system
at this point, then we can write 
$$
\mathbb S_\kappa^2:=\{(x,y,z) \ \! | \ \!  \kappa(x^2+y^2+z^2)+2\kappa^{\frac{1}{2}}z=0\} \ \ {\rm for}\ \ \kappa>0,
$$
$$
\mathbb H_\kappa^2:=\{(x,y,z) \ \! | \ \!  \kappa(x^2+y^2-z^2)+2|\kappa|^{\frac{1}{2}}z=0\hspace{0.2cm}z\geq0\}\ \ {\rm for}\ \ \kappa<0.
$$
Consider now four point masses, $m_i>0,\ i=1,2,3,4$,
whose position vectors, velocities, and accelerations are
given by 
$$
{\bf r}_i=(x_i,y_i,z_i),\ \dot{\bf r}_i=(\dot x_i,\dot y_i,\dot z_i),\ \ddot{\bf r}_i=(\ddot x_i,\ddot y_i,\ddot z_i),\ i=1,2,3,4.
$$ 
Then, as shown in \cite{Diacu04}, the equations of motion take the form
$$
  \begin{cases} 
   \ddot x_i= \sum _{j=1, j\ne i}^N \frac{m_j\Big[{ x}_j-\Big(1-\frac{\kappa r_{ij}^2}{2}\Big){x}_i 
   	\Big]}{\Big(1-\frac{\kappa r_{ij}^2}{4}\Big)^{3/2}r_{ij}^3}- \kappa (\dot{\bf r}_i\cdot \dot{\bf r}_i)x_i  \\ 
   
   \ddot y_i= \sum _{j=1, j\ne i}^N \frac{m_j\Big[{ y}_j-\Big(1-\frac{\kappa r_{ij}^2}{2}\Big){y}_i 
   	\Big]}{\Big(1-\frac{\kappa r_{ij}^2}{4}\Big)^{3/2}r_{ij}^3}- \kappa (\dot{\bf r}_i\cdot \dot{\bf r}_i)y_i \\ 
   \ddot z_i= \sum _{j=1, j\ne i}^N \frac{m_j\Big[{z }_j-\Big(1-\frac{\kappa r_{ij}^2}{2}\Big){z}_i 	\Big]}{\Big(1-\frac{\kappa r_{ij}^2}{4}\Big)^{3/2}r_{ij}^3}- (\dot{\bf r}_i\cdot \dot{\bf r}_i)(\kappa z_i+\sigma|\kappa| ^{1/2}),\ i=1,2,3,4, \end{cases}
$$
where $\sigma=1$ for $\kappa\ge 0$, $\sigma=-1$ for $\kappa<0$, and
$$
r_{ij}:=
\begin{cases}
[(x_i-x_j)^2+(y_i-y_j)^2+(z_i-z_j)^2]^{1/2}\ \
{\rm for}\ \ \kappa\ge0\cr
[(x_i-x_j)^2+(y_i-y_j)^2-(z_i-z_j)^2]^{1/2}\ \
{\rm for}\ \ \kappa<0\cr
\end{cases}
$$
for $i,j\in\{1,2,3,4\}$. The above system has eight constraints, namely
  \begin{equation*}
  \kappa(x_1^2+y_i^2+\sigma z_i^2)+2|\kappa|^{1/2}z_i=0,
  \end{equation*}
  \begin{equation*}
  \kappa{\bf r}_i\cdot \dot{\bf r}_i+|\kappa|^{1/2}\dot z_i=0,\hspace{0.3cm}i=1,2,3,4.
  \end{equation*}
If satisfied at an initial instant, these constraints are satisfied for all time because the sets $\mathbb S_\kappa^2, \R^2$, and $\mathbb H_\kappa^2$ are invariant for the equations of motion, \cite{Diacu03}. Notice that for $\kappa=0$ we recover the classical Newtonian equations of the $4$-body problem on the Euclidean plane, namely 
 \begin{equation*} 
 \ddot{\bf r}_i=\sum_{j=1, j\ne i}^N\frac{m_j({\bf r}_j-{\bf r}_i)}{r_{ij}^3},
 \end{equation*}
 where ${\bf r}_i=(x_i,y_i,0),\ i=1,2,3,4.$  
 
\section{Relative equilibria}
   
Relative equilibria are solutions for which the relative distances remain constant during the motion. We first introduce some coordinates $(\varphi,\omega)$, which were originally used in \cite{Diacu04} for the case $N=3$, to detect relative equilibria on and near the equator of $\mathbb S_\kappa^2$, where $\varphi$ measures the angle from the $x$-axis in the $xy$-plane, while $\omega$  is the height on the vertical $z$-axis. In these new coordinates, the constraints become
 $$x_i^2+y_i^2+\omega_i^2+2\kappa^{-1/2}\omega_i=0,\hspace{0.2cm}i=1,2,3,4.$$
With the notation,
 $$\Omega_i=x_i^2+y_i^2=-\kappa^{-1/2}\omega_i(\kappa^{1/2}\omega_i+2) \geq0, \ \omega_i\in [-2\kappa^{-1/2},0],\ i=1,2,3,4,$$
where equality occurs when the body is at the North or the South Pole of the sphere, the $(\varphi,\omega)$-coordinates are given   by the transformations 
 $$x_i=\Omega_i^{1/2}\cos\varphi_i,\  y_i=\Omega_i^{1/2}\sin\varphi_i.$$
 Thus the equations of motion take the form
$$
 \begin{cases}
 \ddot{\varphi_i}=\Omega_i^{-1/2}\sum _{j=1, j\ne i}^N \frac{ m_j \Omega_j^{1/2}\sin (\varphi_j-\varphi_i)} {\rho_{ij}^3(1-\frac{\kappa \rho_{ij}^2}{4})^{3/2}}-\frac{\dot{\varphi_i}\dot{\Omega_i}}{\Omega_i}\\
  \ddot{\omega_i}=\Omega_i^{-1/2}\sum _{j=1, j\ne i}^N \frac{m_j\Big[ \omega_j+\omega_i+\frac{\kappa \rho_{ij}^2}{2}(\omega_i+\kappa^{-1/2})\Big]} {\rho_{ij}^3 \Big( 1-\frac{\kappa \rho_{ij}^2}{4}\Big)^{3/2}}-(\kappa \omega_i+\kappa^{1/2})(\frac{\dot{\Omega_i^2}}{4\Omega_i}+\dot{\varphi_i}^2\Omega_i+\dot{\omega_i^2}),
 \end{cases}
$$  
where
 $$ 
 \dot{\Omega}_i=-2\kappa^{-1/2}\dot{\omega_i}(\kappa^{1/2}\omega_i+1)
 $$
 $$ 
 \rho_{ij}^2=\Omega_i+\Omega_j-2\Omega_i^{1/2}\Omega_j^{1/2}\cos(\varphi_i-\varphi_j)+(\omega_i-\omega_j)^2,\hspace{0.2cm} i,j=1,2,3,4, \hspace{0.2cm}i\neq j.
$$

\section{Relative equilibria on the equator}

If we restrict the motion of the four bodies to the equator of $\mathbb S_\kappa^2$, then 
$$
\omega_i=-\kappa^{-1/2},\hspace{0.3cm} \dot{\omega_i}=0, \hspace{0.3cm} \Omega_i=\kappa^{-1}, \ i=1,2,3,4,$$ 
and the equations of motion take the form
$$
\ddot{\varphi}_i=\kappa^{3/2}\sum_{j=1,j\ne i}^ 4\frac{m_j\sin(\varphi_j-\varphi_i)}{|\sin(\varphi_j-\varphi_i)|^3}, \ \ i=1,2,3,4.
$$
For the relative equilibria, the angular velocity is the same constant for all masses, so we denote this velocity by $\alpha\neq0$ and take 
$$
 \varphi_1=\alpha t+a_1,\hspace{0.3 cm} \varphi_2=\alpha t+a_2,\hspace{0.3 cm} \varphi_3=\alpha t+a_3,\hspace{0.3 cm} \varphi_4=\alpha t+a_4, 
$$
where $a_1,a_2,a_3,a_4$ are real constants, so
  $$ 
  \ddot{\varphi}_i=0, \ i=1,2,3,4.
  $$
Using the notation
$$
s_1:=\frac {\kappa^{3/2}\sin(\varphi_1-\varphi_2)}{|\sin(\varphi_1-\varphi_2)|^3}, \hspace{0.3cm}  s_2:=\frac{\kappa^{3/2}\sin(\varphi_2-\varphi_3)}{|\sin(\varphi_2-\varphi_3)|^3},\hspace{0.3cm} s_3:=\frac {\kappa^{3/2}\sin(\varphi_3-\varphi_1)}{|\sin(\varphi_3-\varphi_1)|^3},
$$ 	
$$
s_4:=\frac{\kappa^{3/2}\sin(\varphi_4-\varphi_1)}{|\sin(\varphi_4-\varphi_1)|^3},\hspace{0.3cm}  s_5:=\frac {\kappa^{3/2}\sin(\varphi_2-\varphi_4)}{|\sin(\varphi_2-\varphi_4)|^3}, \hspace{0.3cm}   s_6:=\frac{\kappa^{3/2}\sin(\varphi_3-\varphi_4)}{|\sin(\varphi_3-\varphi_4)|^3},
$$
we obtain from the equations of motion that
$$
\begin{cases}
-m_2s_1+m_3s_3+m_4s_4=0\cr
	m_1s_1-m_3s_2-m_4s_5=0\cr
	-m_1s_3+m_2s_2-m_4s_6=0\cr
	-m_1s_4+m_2s_5+m_3s_6=0.
\end{cases}
$$
To have other solutions of the masses than $m_1=m_2=m_3=m_4=0$, the determinant of the
above system must vanish, which is equivalent to
$$
s_1s_6+s_3s_5=s_2s_4.
$$
We have thus proved the following result.
\begin{theorem}
A necessary condition that the quadrilateral inscribed in the equator of $ \mathbb S_\kappa^2$, with the four masses $m_1,m_2,m_3,m_4 >0$ at its vertices, forms a relative equilibrium is that 
$$ 
s_1s_6+s_3s_5=s_2s_4. 
$$
\end{theorem}

\section{Equivalent equations of motion}

Let us now introduce some equivalent equations of motion that
are suitable for the kind of solutions we are seeking. First, rewriting the above constraints as 
\begin{equation*}
\kappa(x_i^2+y_i^2)+(|\kappa|^{1/2}z_i+1)^2=1,
\end{equation*}
and solving explicitly for $z_i$, we obtain
\begin{equation*}
z_i=|\kappa|^{-1/2}[\sqrt{1-\kappa(x_i^2+y_i^2)}-1].
\end{equation*} 
The idea here is to eliminate the four equations involving $z_1,z_2,z_3,z_4$, but they still appear in the terms $r_{ij}^2$ in the form $\sigma(z_i-z_j)^2$ as 
\begin{equation*}
\sigma(z_i-z_j)^2=\frac{\kappa(x_i^2+y_i^2-x_j^2-y_j^2)^2}{\left[\sqrt{1-\kappa(x_i^2+y_i^2)}+\sqrt{1-\kappa(x_j^2+y_j^2)}\right]^2}.
\end{equation*}
The case of physical interest is when $\kappa$ is not far from
zero, so the above expression exist even for $\kappa>0$ under this assumption. Then the equations of motion become
$$ 
\begin{cases} 
\ddot {x_i}= \sum _{j=1, j\ne i}^N \frac{m_j\Big[{ x}_j-\Big(1-\frac{\kappa \rho_{ij}^2}{2}\Big){x}_i 
	\Big]}{\Big(1-\frac{\kappa \rho_{ij}^2}{4}\Big)^{3/2}\rho_{ij}^3}- \kappa (\dot{x_i}^2+\dot{y_i}^2+\kappa B_i)x_i \\ 
\ddot {y_i}= \sum _{j=1, j\ne i}^N \frac{m_j\Big[{ y}_j-\Big(1-\frac{\kappa \rho_{ij}^2}{2}\Big){y}_i 
	\Big]}{\Big(1-\frac{\kappa \rho_{ij}^2}{4}\Big)^{3/2}\rho_{ij}^3}- \kappa (\dot{x_i}^2+\dot{y_i}^2+\kappa B_i)y_i,
\end{cases}
$$
where 
\begin{equation*}
\rho_{ij}^2=(x_i-x_j)^2+(y_i-y_j)^2+\frac{\kappa(A_i-A_j)^2}{(\sqrt{1-\kappa A_i}+\sqrt{1-\kappa A_j})^2},
 \end{equation*}
 \begin{equation*}
 A_i=x_i^2+y_i^2,
 \end{equation*}
\begin{equation*}
B_i=\frac{(x_i\dot{x_i}+y_i\dot{y_i})^2}{1-\kappa(x_i^2+y_i^2)},\hspace{0.4cm}i=1,2,3,4. 
\end{equation*}
It is obvious that for $\kappa=0$ we recover the classical Newtonian equations of motion of the planar 4-body problem.

\section{The case of two negligible masses}

We now consider the case when two out of the four given masses are negligible, $m_3=m_4=0$. Then the equations of motion become 
$$
 \begin{cases} 
\ddot {x}_1= \frac{m_2\Big[{ x}_2-\Big(1-\frac{\kappa \rho_{12}^2}{2}\Big){x}_1 
\Big]}{\Big(1-\frac{\kappa \rho_{12}^2}{2}\Big)^{3/2}\rho_{12}^3}- \kappa (\dot{x_1}^2+\dot{y_1}^2+\kappa B_1)x_1 \\ 

\ddot {y}_1=  \frac{m_2\Big[{ y}_2-\Big(1-\frac{\kappa \rho_{12}^2}{2}\Big){y}_1 
	\Big]}{\Big(1-\frac{\kappa \rho_{12}^2}{4}\Big)^{3/2}\rho_{12}^3}- \kappa (\dot{x_1}^2+\dot{y_1}^2+\kappa B_1)y_1 \\
 \end{cases}
$$
$$
\begin{cases}
\ddot {x}_2= \frac{m_1\Big[{ x}_1-\Big(1-\frac{\kappa \rho_{12}^2}{2}\Big){x}_2 
	\Big]}{\Big(1-\frac{\kappa \rho_{12}^2}{4}\Big)^{3/2}\rho_{12}^3}- \kappa (\dot{x_2}^2+\dot{y_2}^2+\kappa B_2)x_2 \\ 

\ddot {y}_2=  \frac{m_1\Big[{ y}_1-\Big(1-\frac{\kappa \rho_{12}^2}{2}\Big){y}_2 
	\Big]}{\Big(1-\frac{\kappa \rho_{12}^2}{4}\Big)^{3/2}\rho_{12}^3}- \kappa (\dot{x_2}^2+\dot{y_2}^2+\kappa B_2)y_2 \\
\end{cases}
$$
$$
\begin{cases}
\ddot {x}_3= \frac{m_1\Big[{ x}_1-\Big(1-\frac{\kappa \rho_{13}^2}{2}\Big){x}_3 
  	\Big]}{\Big(1-\frac{\kappa \rho_{13}^2}{4}\Big)^{3/2}\rho_{13}^3}+ \frac{m_2\Big[{ x}_2-\Big(1-\frac{\kappa \rho_{32}^2}{2}\Big){x}_3 
  	\Big]}{\Big(1-\frac{\kappa \rho_{23}^2}{4}\Big)^{3/2}\rho_{23}^3}- \kappa (\dot{x_3}^2+\dot{y_3}^2+\kappa B_3)x_3 \\ 
  
 \ddot {y}_3= \frac{m_1\Big[{ y}_1-\Big(1-\frac{\kappa \rho_{13}^2}{2}\Big){y}_3 
  	\Big]}{\Big(1-\frac{\kappa \rho_{13}^2}{4}\Big)^{3/2}\rho_{13}^3}+ \frac{m_2\Big[{ y}_2-\Big(1-\frac{\kappa \rho_{32}^2}{2}\Big){y}_3 
  	\Big]}{\Big(1-\frac{\kappa \rho_{32}^2}{2}\Big)^{3/2}\rho_{32}^3}- \kappa (\dot{x_3}^2+\dot{y_3}^2+\kappa B_3)y_3 \\
\end{cases}
$$
$$
  \begin{cases}
 \ddot {x}_4= \frac{m_1\Big[{ x}_1-\Big(1-\frac{\kappa \rho_{14}^2}{2}\Big){x}_4 
 	\Big]}{\Big(1-\frac{\kappa \rho_{14}^2}{4}\Big)^{3/2}\rho_{14}^3}+ \frac{m_4\Big[{ x}_4-\Big(1-\frac{\kappa \rho_{42}^2}{2}\Big){x}_2 
 	\Big]}{\Big(1-\frac{\kappa \rho_{42}^2}{4}\Big)^{3/2}\rho_{42}^3}- \kappa (\dot{x_4}^2+\dot{y_4}^2+\kappa B_4)x_4 \\ 
 \ddot {y}_4= \frac{m_1\Big[{ y}_1-\Big(1-\frac{\kappa \rho_{14}^2}{2}\Big){y}_4 
 	\Big]}{\Big(1-\frac{\kappa \rho_{14}^2}{4}\Big)^{3/2}\rho_{14}^3}+ \frac{m_4\Big[{ y}_4-\Big(1-\frac{\kappa \rho_{42}^2}{2}\Big){y}_2 
 	\Big]}{\Big(1-\frac{\kappa \rho_{42}^2}{4}\Big)^{3/2}\rho_{42}^3}- \kappa (\dot{x_4}^2+\dot{y_4}^2+\kappa B_4)y_4, 
\end{cases}
$$
where
 $ \rho_{ij}^2=\rho_{ji}^2, \ i\neq j$,
$$
\rho_{ij}^2=(x_i-x_j)^2+(y_i-y_j)^2+\frac{\kappa(x_i^2+y_i^2-x_j^2-y_j^2)^2}{[\sqrt{1-\kappa(x_i^2+y_i^2)}+\sqrt{1-\kappa(x_j^2+y_j^2)}]^2}. 
$$
We can now show that when $m_1=m_2=:m>0$ and $m_3=m_4=0$, square-like relative equilibria, i.e.\ equilateral equiangular quadrilaterals, always exist on $\mathbb S_\kappa^2$, but not on $\mathbb H_\kappa^2$.
\begin{theorem}
In the curved 4-body problem,
assume that $m_1=m_2=:m>0$ and $m_3=m_4=0$. Then, in $\mathbb S_\kappa^2$, there is a circle of radius $r<\kappa^{-1/2}$, parallel with the $xy$-plane, such that a square configuration inscribed in this circle, with $m_1, m_2$ at the opposite ends of one diagonal and $m_3, m_4$ at the opposite ends of the other diagonal, forms a relative equilibrium. But in $\mathbb H_\kappa^2$, there is no such solution.
\end{theorem}

	\begin{figure}
		\begin{tikzpicture}[scale=0.7]
		\tikzstyle{every node}=[font=\large]
		\shade[ball color = gray!30, opacity = 0.3] (0,0) circle (4cm);
		\draw (0,0) circle (4cm);
		\draw (-4,0) arc (180:360:4 and 0.6);
		\draw[dashed] (4,0) arc (0:180:4 and 0.6);
		\draw (-2.8,2.8) arc (180:360:2.8 and 0.5);
		\draw[dashed] (2.8,2.8) arc (0:180:2.8 and 0.5);
		\fill[fill=black] (0,0) circle (1pt);
	\fill (-2.8,2.8) circle (8pt)node[above left]{$m_1$};
	\fill (2.8,2.8) circle (8pt)node[above right]{$m_2$};
		\fill (0,2.3) circle (5pt)node[above]{$m_4$};
		\fill (0,3.3) circle (5pt)node[above]{$m_3$};		
		\end{tikzpicture}
		\caption{ The case of  2 equal masses and 2  negligible masses.}
	\end{figure}
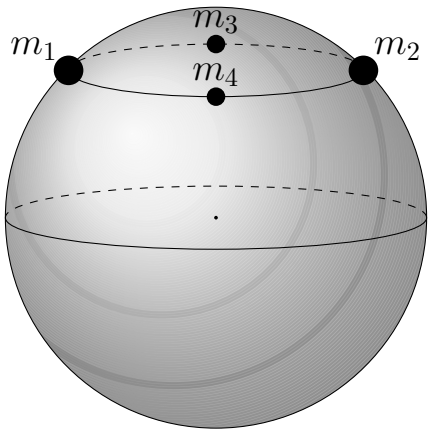

\begin{proof}
We must check the existence of a solution of the form 
$$ {\bf q }=(q_1,q_2,q_3,q_4)\in \mathbb S_\kappa^2,\hspace{0.3cm} {\bf q_i} =(x_i,y_i),\hspace {0.2cm}i=1,2,3,4.$$
\begin{equation*}
  x_1=r\cos\alpha t,\hspace{0.5cm} y_1=r \sin\alpha t, \hspace{0.9cm}  
  \end{equation*}   
 \begin{equation*}
 x_2=-r\cos\alpha t,\hspace{0.5cm} y_2=-r \sin\alpha t, \hspace{0.5cm} 
 \end{equation*}
  \begin{equation*}
  x_3= r\cos(\alpha t+\pi/2) =-r\sin\alpha t,\hspace{0.5cm} y_3= r\sin(\alpha t+\pi/2)=r \cos\alpha t, \hspace{0.5cm}  
 \end{equation*} 
 \begin{equation*}
 x_4= -r\cos(\alpha t+\pi/2) =r\sin\alpha t,\hspace{0.8cm} y_4=- r\sin(\alpha t+\pi/2)=-r \cos\alpha t, \hspace{0.5cm}  
 \end{equation*}
where 
$$
\hspace{0.5cm} x_i^2+y_i^2=r^2,\ \
\rho^2 = \rho_{13}^2=\rho_{14}^2=\rho_{23}^2=\rho_{24}^2= 2r^2, \ \  
\rho_{12}^2=\rho_{34}^2=4r^2. 
$$
Substituting these expressions into the system, the first four equations lead us to
$$
\alpha^2=\frac{m}{4r^3(1-\kappa r^2)^{3/2}},
$$
whereas the last four equations yield
$$
\alpha^2=\frac{2m(1-\frac{\kappa \rho^2}{2})}{\rho^3(1- \frac{\kappa \rho^2}{4})^{3/2}(1-\kappa r^2)}.
$$ 
So, to have a solution, the equation 
$$
 \frac{m}{4r^3(1-\kappa r^2)^{3/2}}=\frac{2m(1-\frac{\kappa\rho^2}{2})}{\rho^3(1- \frac{\kappa \rho^2}{4})^{3/2}(1-\kappa r^2 )}
$$
must be satisfied. This equation is equivalent to 
$$
\frac{1}{8r^3(1-\kappa r^2)^{3/2}}=\frac{1}{2\sqrt{2} r^3(1- \frac{\kappa r^2}{2})^{3/2}},
$$
which leads to 
$$
3\kappa r^2=2.
$$
Obviously, in the case of $\mathbb H_\kappa^2$, we have 
$\kappa<0$, so this equation has no solutions. For $\mathbb S_\kappa^2$, it leads to
$$
r=\sqrt{2/3}\kappa^{-1/2}.
$$
Since $r<\kappa^{-1/2}$, such a solution always exists in
$\mathbb S_\kappa^2$.     
\end{proof}

	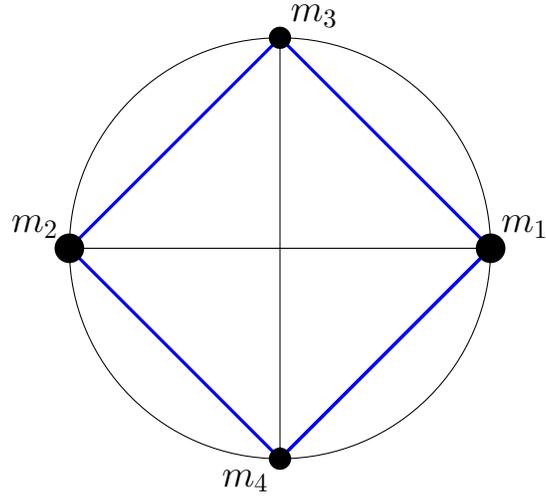
\begin{figure}
		\begin{tikzpicture}[scale=0.7]
		\tikzstyle{every node}=[font=\large]
		\draw (0, 0) circle (4);	
		\draw [->](0,0)--(90:4);
		\draw [->](0,0)--(0:4);
		\draw [->] [draw=blue, very thick](90:4)--(0:4);
		\draw [->](0,0)--(0:-4);
		\draw [->][draw=blue, very thick](90:4)--(0:-4);
		\draw[->] [draw=blue, very thick](0:4)--(270:4);
		\draw [->](0,0)--(270:4);
		\draw [draw=blue, very thick] [->](0:-4)--(270:4);
		
		\draw [draw=blue, very thick] [->](0:4)--(270:4);
		
\fill (0:4) circle (8pt)node[above right]{$m_1$};
\fill  (90:4) circle (6pt)node[above right]{$m_3$};
\fill  (0:-4) circle (8pt)node[above left]{$m_2$};
\fill  (270:4) circle (6pt)node[below left]{$m_4$};
		\end{tikzpicture}
		\caption{The case of two equal masses and two negligible masses.}
	\end{figure}

\section{The case of one negligible mass}

Let $m_1,m_2,m_3 >0$ and assume that $m_4=0$. Then the equations of motion take the form
$$
\begin{cases}
\ddot {x_1}= \frac{m_2\Big[{ x}_2-\Big(1-\frac{\kappa \rho_{12}^2}{2}\Big){x}_1 
	\Big]}{\Big(1-\frac{\kappa \rho_{12}^2}{4}\Big)^{3/2}\rho_{12}^3}+ \frac{m_3\Big[{ x}_3-\Big(1-\frac{\kappa \rho_{31}^2}{2}\Big){x}_1 
	\Big]}{\Big(1-\frac{\kappa \rho_{31}^2}{4}\Big)^{3/2}\rho_{31}^3}- \kappa (\dot{x_1}^2+\dot{y_1}^2+\kappa B_1)x_1 \\ 

\ddot {y_1}= \frac{m_2\Big[{ y}_2-\Big(1-\frac{\kappa \rho_{12}^2}{2}\Big){y}_1 
	\Big]}{\Big(1-\frac{\kappa \rho_{12}^2}{4}\Big)^{3/2}\rho_{12}^3}+ \frac{m_3\Big[{ y}_3-\Big(1-\frac{\kappa \rho_{31}^2}{2}\Big){y}_1 
	\Big]}{\Big(1-\frac{\kappa \rho_{31}^2}{4}\Big)^{3/2}\rho_{31}^3}- \kappa (\dot{x_1}^2+\dot{y_1}^2+\kappa B_1)y_1 \\ 
\end{cases} 
$$
$$
\begin{cases}
\ddot {x_2}= \frac{m_1\Big[{ x}_1-\Big(1-\frac{\kappa \rho_{12}^2}{2}\Big){x}_2 
	\Big]}{\Big(1-\frac{\kappa \rho_{12}^2}{4}\Big)^{3/2}\rho_{12}^3}+ \frac{m_3\Big[{ x}_3-\Big(1-\frac{\kappa \rho_{32}^2}{2}\Big){x}_2 
	\Big]}{\Big(1-\frac{\kappa \rho_{32}^2}{4}\Big)^{3/2}\rho_{32}^3}- \kappa (\dot{x_2}^2+\dot{y_2}^2+\kappa B_2)x_2 \\ 
 \ddot {y_2}= \frac{m_1\Big[{ y}_1-\Big(1-\frac{\kappa \rho_{12}^2}{2}\Big){y}_2 
	\Big]}{\Big(1-\frac{\kappa \rho_{12}^2}{4}\Big)^{3/2}\rho_{12}^3}+ \frac{m_3\Big[{ y}_3-\Big(1-\frac{\kappa \rho_{32}^2}{2}\Big){y}_2 
	\Big]}{\Big(1-\frac{\kappa \rho_{32}^2}{4}\Big)^{3/2}\rho_{32}^3}- \kappa (\dot{x_2}^2+\dot{y_2}^2+\kappa B_2)y_2 
\end{cases} 
$$
$$
\begin{cases}
\ddot {x_3}= \frac{m_1\Big[{ x}_1-\Big(1-\frac{\kappa \rho_{13}^2}{2}\Big){x}_3 
	\Big]}{\Big(1-\frac{\kappa \rho_{13}^2}{4}\Big)^{3/2}\rho_{13}^3}+ \frac{m_2\Big[{ x}_2-\Big(1-\frac{\kappa \rho_{32}^2}{2}\Big){x}_3 
	\Big]}{\Big(1-\frac{\kappa \rho_{32}^2}{4}\Big)^{3/2}\rho_{32}^3}- \kappa (\dot{x_3}^2+\dot{y_3}^2+\kappa B_3)x_3 \\ 
\ddot {y_3}= \frac{m_1\Big[{ y}_1-\Big(1-\frac{\kappa \rho_{13}^2}{2}\Big){y}_3 
	\Big]}{\Big(1-\frac{\kappa \rho_{13}^2}{4}\Big)^{3/2}\rho_{13}^3}+ \frac{m_2\Big[{ y}_2-\Big(1-\frac{\kappa \rho_{32}^2}{2}\Big){y}_3 
	\Big]}{\Big(1-\frac{\kappa \rho_{32}^2}{4}\Big)^{3/2}\rho_{32}^3}- \kappa (\dot{x_3}^2+\dot{y_3}^2+\kappa B_3)y_3 
\end{cases} 
$$
$$
\begin{cases}
\ddot {x_4}= \frac{m_1\Big[{ x}_1-\Big(1-\frac{\kappa \rho_{14}^2}{2}\Big){x}_4 
	\Big]}{\Big(1-\frac{\kappa \rho_{14}^2}{4}\Big)^{3/2}\rho_{14}^3}+ \frac{m_2\Big[{ x}_2-\Big(1-\frac{\kappa \rho_{42}^2}{2}\Big){x}_4 
	\Big]}{\Big(1-\frac{\kappa \rho_{42}^2}{4}\Big)^{3/2}\rho_{42}^3}+ \frac{m_3\Big[{ x}_3-\Big(1-\frac{\kappa \rho_{43}^2}{2}\Big){x}_4 
	\Big]}{\Big(1-\frac{\kappa \rho_{43}^2}{4}\Big)^{3/2}\rho_{43}^3}\\
	\hfill -\kappa (\dot{x_4}^2+\dot{y_4}^2+\kappa B_4)x_4 \\ 
\ddot {y_4}= \frac{m_1\Big[{ y}_1-\Big(1-\frac{\kappa \rho_{14}^2}{2}\Big){y}_4 
	\Big]}{\Big(1-\frac{\kappa \rho_{14}^2}{4}\Big)^{3/2}\rho_{14}^3}+ \frac{m_2\Big[{ y}_2-\Big(1-\frac{\kappa \rho_{42}^2}{2}\Big){y}_4 
	\Big]}{\Big(1-\frac{\kappa \rho_{42}^2}{4}\Big)^{3/2}\rho_{42}^3}+ \frac{m_3\Big[{ y}_3-\Big(1-\frac{\kappa \rho_{43}^2}{2}\Big){y}_4 
	\Big]}{\Big(1-\frac{\kappa \rho_{43}^2}{4}\Big)^{3/2}\rho_{43}^3}\\
	\hfill - \kappa (\dot{x_4}^2+\dot{y_4}^2+\kappa B_4)y_4.  
\end{cases} 
$$
We will next show that if the non-negligible masses are equal, then there exist some kite-shaped relative equilibria. 
\begin{theorem}
Consider the curved 4-body problem with masses $m_1=m_2=m_3:=m>0$ and $m_4=0$. Then, in $\mathbb S_\kappa^2$, there exists at least one kite-shaped relative equilibrium for which the equal masses lie at the vertices of an equilateral triangle, whereas the negligible mass is at the intersection of the extension of one height of the triangle with the circle on which all the bodies move. In $\mathbb H_\kappa^2$, however, there are no such kite-shaped relative equilibria.
\end{theorem}

	  \begin{figure}
	  	\begin{tikzpicture}[scale=0.7]
		\tikzstyle{every node}=[font=\large]
	  	\draw (0, 0) circle (4);	  	
	  	\draw [->](0,0)--(90:4);
	  	\draw [->](0,0)--(-30:4);
	\draw [->] [draw=blue, very thick](90:4)--(-30:4);
	\draw [->](0,0)--(210:4);
	\draw [->][draw=blue, very thick](90:4)--(210:4);
	\draw[->] [draw=blue, very thick](210:4)--(-30:4);
	\draw [->](0,0)--(270:4);
	\draw [draw=blue, very thick] [->](-30:4)--(270:4);
	\draw [draw=blue, very thick] [->](210:4)--(270:4);	  	
\fill[black] (-30:4) circle (8pt)node[below right]{$m_1$};
\fill [black] (90:4) circle (8pt)node[above right]{$m_2$};
\fill [black] (210:4) circle (8pt)node[below left]{$m_3$};
\fill [black] (270:4) circle (6pt)node[below left]{$m_4$};
	  	\end{tikzpicture}
	  	\caption{A kite configuration of 3 equal masses and one negligible mass.}
	  \end{figure}
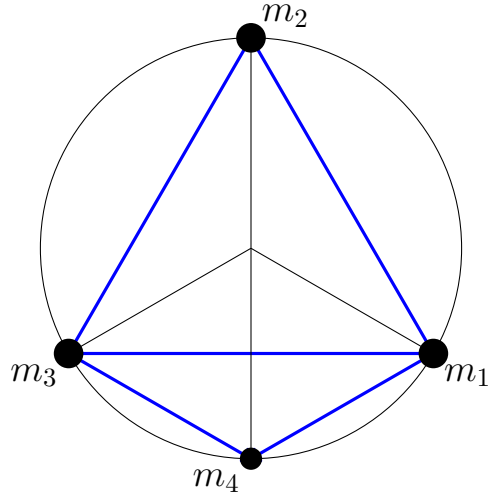

\begin{proof}
We will check a solution of the form
 \begin{equation*}
  x_1=r\cos\alpha t,\hspace{0.5cm} y_1=r \sin\alpha t, \hspace{1cm}  
  \end{equation*}   
  \begin{equation*}
  x_2=r\cos\Big (\alpha t+\frac{2\pi}{3}\Big),\hspace{0.5cm} y_2=r \sin \Big (\alpha t+\frac{2\pi}{3}\Big) 
  \end{equation*}
  \begin{equation*}
  x_3= r \cos\Big (\alpha t+\frac{4\pi}{3}\Big),\hspace{0.5cm} y_3= r \sin\Big (\alpha t+\frac{4\pi}{3}\Big), \hspace{0.5cm} 
  \end{equation*} 
  \begin{equation*}
  x_4= r \cos\Big (\alpha t-\frac{\pi}{3}\Big),\hspace{0.5cm} y_4= r \sin\Big (\alpha t-\frac{\pi}{3}\Big), \hspace{0.5cm} 
  \end{equation*}
where 
$$ 
\rho_{12}^2= \rho_{13}^2=\rho_{23}^2=3 r^2,\ \
\rho_{43}^2=\rho_{41}^2= r^2,\ \
\rho_{24}^2=4r^2.
$$
Substituting these expressions into the above system, we 
are led to the conclusion that the following two equations
must be satisfied,
$$
\alpha^2=\frac{m}{\sqrt{3} r^3 (1-\frac{3\kappa r^2}{4})^{3/2}},
$$   
$$
\alpha^2=\frac{m}{4 r^3 (1-\kappa r^2)^{3/2}}+\frac{m}{ r^3 (1-\frac{\kappa r^2}{4})^{3/2}}.
$$
Comparing these equations we obtain the condition for the existence of the kite-shaped relative equilibria, 
$$
\frac{1}{\sqrt{3} (1-\frac{3\kappa r^2}{4})^{3/2}}=\frac{1}{4  (1-\kappa r^2)^{3/2}}+\frac{1}{ (1-\frac{\kappa r^2}{4})^{3/2}}.
$$
Straightforward computations show that $r$ is a solution of this equation if it is a root of the polynomial 
$$
P(r)=a_{24}r^{24}+a_{22}r^{22}+a_{20}r^{20}+a_{18}r^{18}+a_{16}r^{16}+a_{14}r^{14}+a_{12}r^{12}+$$
$$ a_{10}r^{10}+a_{8}r^{8}+a_{6}r^{6}+a_4r^4+a_2r^2+a_0,
$$
$$
a_{24}=\frac{6697290145}{16777216}\kappa^{12},\ \ a_{22}=-\frac{2884257825}{524288}\kappa^{11},
\ \ a_{20}=\frac{18063189465}{524288}\kappa^{10},
$$
$$
a_{18}=-\frac{4241985935}{32768}\kappa^9, \ \ \ a_{16}=\frac{21267471735}{65536}\kappa^8,
$$
$$ 
a_{14}=-\frac{584429805}{1024}\kappa^7,\ \ \ a_{12}=\frac{737853351}{1024}\kappa^6,\ \ \ 
a_{10}=-\frac{41995431}{64}\kappa^5,$$
$$
 a_{8}=\frac{109080063}{256}\kappa^4,\ \ \ a_6=-\frac{1530101}{8} \kappa^3,
$$
$$
a_{4}=\frac{446217}{8}\kappa^2, \ \ \ a_{2}=-9318 \kappa,\ \ \ a_0=649$$
that belongs to the interval $r\in (0,\kappa^{-1/2})$ for $\mathbb S_\kappa^2$, but needs only to be positive for $\mathbb H_\kappa^2$. To find out if we have such a root, we make the substitution $ x=r^2$, and obtain the polynomial
$$
Q(x)=a_{24}x^{12}+a_{22}x^{11}+a_{20}x^{10}+a_{18}x^{9}+a_{16}x^{8}+a_{14}x^{7}+a_{12}x^{6}+$$
$$ a_{10}x^{5}+a_{8}x^{4}+a_{6}x^{3}+a_4x^2+a_2x+a_0.
$$
By Descartes's rule of signs the number of positive roots depends on the number of changes of sign of the coefficients, which in turn depends on the sign of $\kappa$. So let us discuss the two cases separately.

In $\mathbb S_\kappa^2$, i.e.\ for $\kappa>0$, there are twelve changes of sign, so $Q$ can have twelve, ten, eight, six, four, two, or zero positive roots, so this does not guarantee the existence of a positive root. However, we can notice that $Q( \frac{\kappa^{-1}}{2}) = -2.4959<0$ and $Q(0)=649>0$, so a root must exist for $x\in(0,\kappa^{-1/2})$, i.e.\ for $r\in(0,\kappa^{-1})$, a remark that proves the existence of at least one kite-shaped relative equilibrium.


In $\mathbb H_\kappa^2$, i.e.\ for $\kappa<0$, we seek a positive root of $Q$. But for $\kappa<0$, there is no sign change in $Q(x)$, so the polynomial has no positive roots. Therefore there are no kite solutions in $\mathbb H_\kappa^2$. This remark completes the proof.
\end{proof}

\bigskip
\noindent{\bf Acknowledgment.} Florin Diacu did most of the work on this paper while visiting the Yangtze Center of Mathematics at Sichuan University as Distinguished Foreign Professor in April-May 2017. He was also supported in part by a grant from the Yale-NUS College at the National University of Singapore and an NSERC of Canada Discovery Grant. Sawsan Alhowaity was funded by a scholarship from the University of Shaqra, Saudia Arabia, towards the completion of her doctoral degree at the University of Victoria in Canada.


\end{document}